\documentclass[preprintnumbers, floatfix, preprintnumbers, letterpaper, twocolumn, superscriptaddress,nofootinbib]{revtex4}
\usepackage{graphicx}
\usepackage{microtype}
\usepackage{amsmath}
\usepackage{amssymb}
\usepackage{subfigure}
\usepackage{hyperref}
\usepackage{url}
\usepackage{xcolor}
\usepackage{color}
\usepackage{mathrsfs}
\usepackage{calrsfs}
\usepackage{amsfonts}
\usepackage{tabularx}
\usepackage{eucal}
\usepackage{latexsym}
\usepackage{ragged2e}
\usepackage{epsfig}
\usepackage{textcomp}
\usepackage{float}

\usepackage{caption}
\DeclareCaptionJustification{justified}{\leftskip=0pt \rightskip=0pt \parfillskip=0pt plus 1fil}
\definecolor{vividviolet}{rgb}{0.62, 0.0, 1.0}
\definecolor{amaranth}{rgb}{0.9, 0.17, 0.31}
\definecolor{palatinateblue}{rgb}{0.15, 0.23, 0.89}
\definecolor{brightpink}{rgb}{1.0, 0.0, 0.5}
\definecolor{cornflowerblue}{rgb}{0.39, 0.58, 0.93}
\definecolor{deepcarminepink}{rgb}{0.94, 0.19, 0.22}
\definecolor{radicalred}{rgb}{1.0, 0.21, 0.37}

\hypersetup{ linktoc=all,
	colorlinks, linkcolor={palatinateblue},
	citecolor={brightpink}, urlcolor={amaranth}
}

\graphicspath{{Images/}}

\renewcommand{\d}[1]{\ensuremath{\operatorname{d}\!{#1}}}

\def\sideremark#1{\ifvmode\leavevmode\fi\vadjust{\vbox to0pt{\vss
			\hbox to 0pt{\hskip\hsize\hskip1em
				\vbox{\hsize1.3cm\tiny\raggedright\pretolerance10000
					\noindent #1\hfill}\hss}\vbox to8pt{\vfil}\vss}}}%
\def\beq{\begin{equation}}
\def\eeq{\end{equation}}

\begin{document}
\title{Lotka-Volterra Models for Extraterrestrial Self-Replicating Probes}

\author{Yifan Chen}
	\affiliation{Center for Gravitation and Cosmology, College of Physical Science and Technology, Yangzhou University, \\180 Siwangting Road, Yangzhou City, Jiangsu Province  225002, China}
	
\author{Jiayi Ni}
\affiliation{Center for Gravitation and Cosmology, College of Physical Science and Technology, Yangzhou University, \\180 Siwangting Road, Yangzhou City, Jiangsu Province  225002, China}

\author{Yen Chin \surname{Ong}}
\email{ycong@yzu.edu.cn}
\affiliation{Center for Gravitation and Cosmology, College of Physical Science and Technology, Yangzhou University, \\180 Siwangting Road, Yangzhou City, Jiangsu Province  225002, China}
\affiliation{Shanghai Frontier Science Center for Gravitational Wave Detection, School of Aeronautics and Astronautics, Shanghai Jiao Tong University, Shanghai 200240, China}

\begin{abstract}
A sufficiently advanced extraterrestrial civilization can send out a swarm of self-replicating probes for space exploration. Given the fast-growing number of such a probe, even if there is only one extraterrestrial civilization sending out such probes in the Milky Way galaxy, we should still expect to see them. The fact that we do not consists part of the Fermi paradox. The suggestion that self-replicating probes will eventually mutate to consume their progenitors and therefore significantly reduce the number of total probes has been investigated and dismissed in the literature. In this work, we re-visit this question with a more realistic Lotka-Volterra model, and show that mutated probes would drive the progenitor probes into ``extinction'', thereby replacing them to spread throughout the galaxy. Thus, the efficiency of mutated probes in reducing the total number of self-replicating probes is even less than previously thought. As part of the analysis, we also suggest that, somewhat counter-intuitively, in designing self-replicating probes, one should not program them to stop replicating when sufficient mutation causes the probes to fail to recognize the progenitor probes as ``self''. 
\end{abstract} 

\maketitle

\section{Introduction: Fermi Paradox and Self-Replicating Probes}

While it is still unclear how widespread intelligent life is across the Universe, a recent optimistic Monte Carlo estimate gives about 42000 communicating extraterrestrial intelligent civilizations in our Milky Way galaxy \cite{2204.05479} alone. If there are so many extraterrestrial civilizations, why haven't we found any of them? This is the famous Fermi paradox \cite{book, 0907.3432, 1104.0624}, or ``the Great Silence'' \cite{silence}. There are of course no shortage of possible explanations (assuming the emergence of life is not rare \cite{rare}), such as: they may not keen on communicating (perhaps for fear of attracting unwanted attention from hostile species \cite{darkforest}); civilizations generally tend to annihilate themselves in one way or another \cite{1507.08530}; we are not considered technologically advanced or interesting enough to be communicated with and are deliberately isolated \cite{zoo}, just to name a few. All the different scenarios we can think of concerning the Fermi paradox is of course necessarily anthropocentric -- we do not know how an alien species might think. 

However, the lack of any confirmed signal from deep space is not the only puzzling aspect of the mystery. We expect that any sufficiently advanced civilization would be able to design and manufacture self-replicating probes \cite{sagan, AS, 1111.6131,1206.0953,1909.05078,dyson} (SRPs, also known as ``replicators'' or von Neumann probes \cite{von}), perhaps even ``nanobots'' that can be launced {en masse} into interstellar space to explore the Universe.  It has been argued that small-scale, partially self-replicating probes are feasible near-term even to us \cite{2005.12303}.
In addition to pure exploration, SRPs can be used for conquest -- obliterating every other civilizations (viewed as competitors or potential threats) and taking over their resources \cite{berserk}. 
SRPs make the Fermi paradox a lot worse, since the fact that such a replication process (assuming virtually unlimited resources) would be exponential, means that even if there is only \emph{one} single civilization that is doing this (conservative estimate gives $O(1)$ to $O(10^2)$ advanced civilizations in the Milky Way \cite{Wallenhorst, 2004.03968}), their probes should have spread throughout the galaxy in a reasonable amount of time (between $O(10^6)$ to $O(10^9)$ years) \cite{1111.6131,Jones}. Given how young we are on the galactic stage, we should have been surrounded by SRPs. The fact that we do not requires some explanations.

One possibility is that SRPs would eventually ``mutate'' \cite{CH, forgan}. While minor mistakes in the replicating process would only create sub-variant probes, a bigger error might lead to a new population of SRPs that can no longer recognize the progenitor probe as ``self'' (as opposed to being ``foreign''), and instead treat it as resource to be consumed. We will refer to such mutated SRPs as ``predators'', while the rest as ``preys''. This comparison with biological system is not entirely an analogy, since SRPs could be partially biological in the sense that they are artificial proteomic or molecular machines \cite{drexler, ellery}. The idea is that predation might reduce the total number of SRPs and thus explains why we have not detected them. In \cite{forgan}, the author investigated this scenario by modeling the spread of the SRPs from one star system to another with minimal spanning trees, and the dynamics between the predator and prey SRPs with the Lotka-Volterra model. The conclusion reached was that this scenario cannot sufficiently reduce the prey SRPs, and therefore the Fermi paradox remains unsolved. 

As we will explain in the next section, there are some shortcomings in the assumptions in \cite{forgan}, and the modeling of predator-prey system with Lotka-Volterra equations require some care in what kind of terms should be included. We then propose a more realistic equation to study the species interaction. Our model makes it clear why predator SRPs cannot resolve the Fermi paradox.
\newpage
\section{Lotka-Volterra Models for Replicators}

Aside from the migration effect from one star system to another, the Lotka-Volterra model assumed in \cite{forgan} is essentially the standard ``plain vanilla'' form\footnote{A finite capacity term, $K >0$, can be included, as was done in \cite{forgan}, so that $\alpha x$ becomes $\alpha x(1-x/K)$, but since $K$ is assumed to be very large anyway, we need not include it in our analysis.}:
\begin{equation}\label{vanilla}
\left\{ 
\begin{aligned}
\frac{\d x}{\d t} &= \alpha x - \beta xy, \notag \\
\frac{\d y}{\d t} &= \delta xy - \gamma y,
\end{aligned}\right.
\end{equation}
where $x=x(t), y=y(t)$ denote the number of preys and predators, respectively. The coefficients $\alpha,\beta,\delta, \gamma$ are assumed to be positive. They govern the growth rate of preys, the reduction in the prey population due to predation, the increase in the predator population due to predation, and the natural death rate of the predator species, respectively. One may ask why the natural death rate of the prey species is not included. The reason is straightforward: we could add such a term $-\alpha'x$ (where $\alpha'>0$), but in a healthy population, the growth rate should be larger than the death rate anyway, so the overall effect is $(\alpha-\alpha')x$, where $\alpha - \alpha' > 0$. So it is easier to absorb $\alpha'$ and just deal with $\alpha$. As the coefficients are unknown, we are only interested in the qualitative features. The units for time $t$ is therefore left unspecified.

\begin{figure}
\begin{center}
\includegraphics[width=0.40\textwidth]{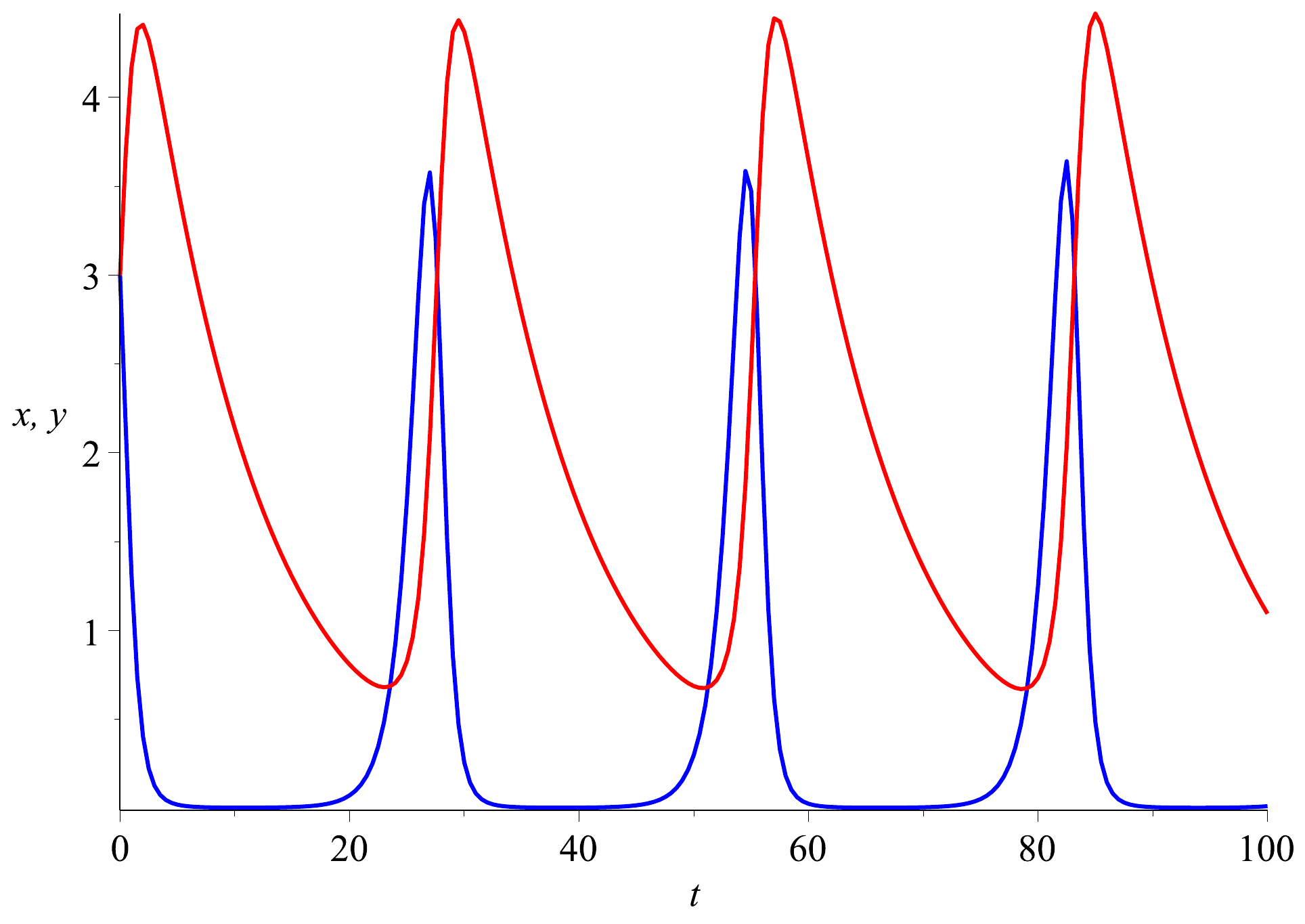}
\caption{The oscillatory number of the predator ($y$, the upper curve in red) and prey ($x$, the lower curve in blue) populations in the ``plain-vanilla'' predator-prey model. Here we use $\alpha=1,\beta=0.5,\gamma=0.2,\delta=0.1$. The initial conditions are $x(0)=y(0)=3$.}
\label{fig1}
\end{center}
\end{figure}

\begin{figure}
\begin{center}
\includegraphics[width=0.42\textwidth]{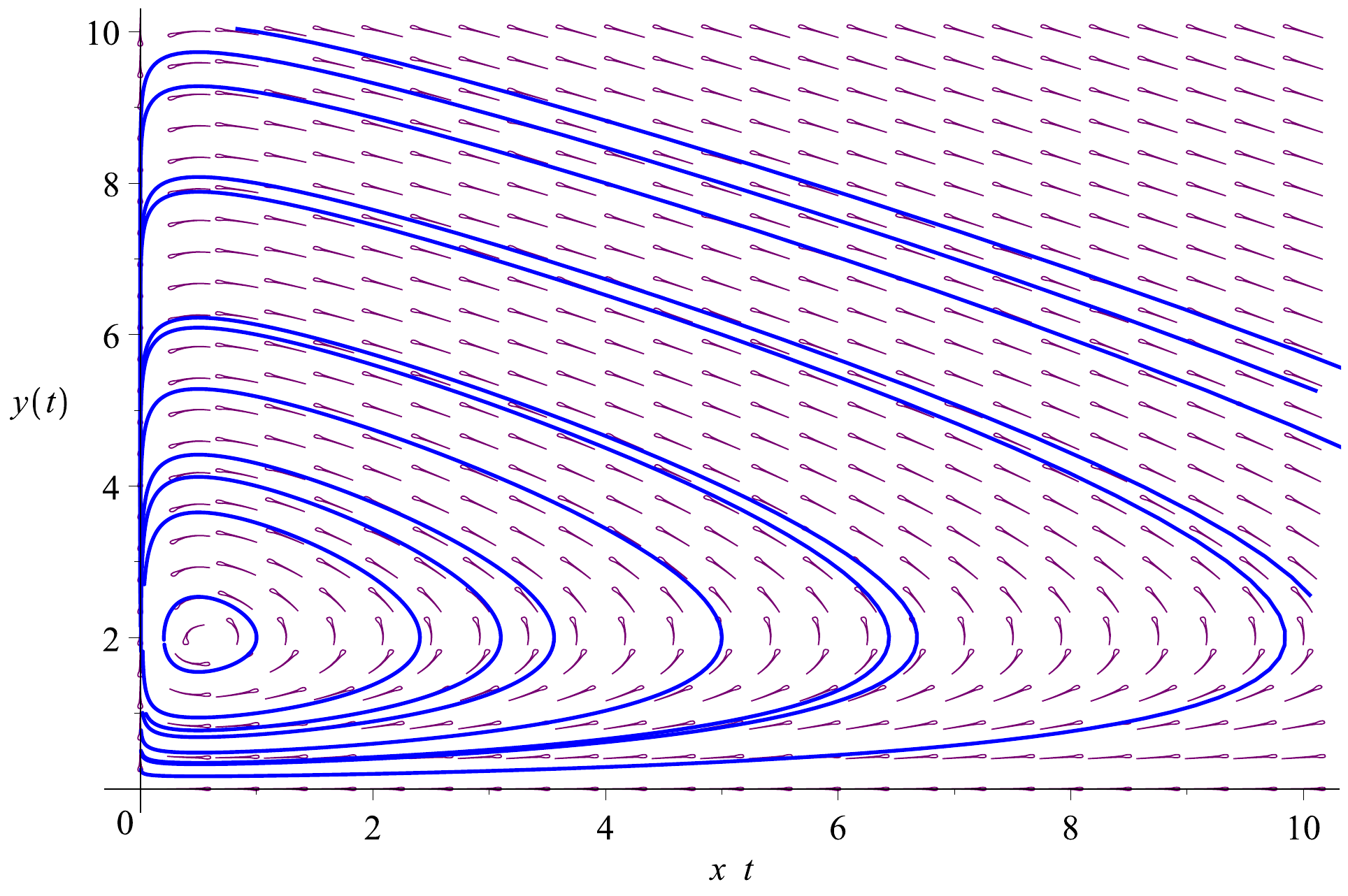}
\caption{The phase diagram of the ``plain-vanilla'' predator-prey model with $\alpha=1,\beta=0.5,\gamma=0.2,\delta=0.1$. There exists a center surrounded by closed orbits.}
\label{fig2}
\end{center}
\end{figure}

The behavior of this ODE system is well-known: the two species are locked into oscillatory cycles, an example is shown in Fig.(\ref{fig1}). As the number of predators increase, the prey population dwindles, and as the preys are reduced, the predators run out of food and they themselves start to die off. The prey population then recovers and the cycle renews. In the associated phase diagram, the only fixed point is a center (Fig.(\ref{fig2})). While this model is rather simplified and suffers from some problems (such as the ``attofox'' problem \cite{atto} -- the number of the species can get very small, which in the real world could have led to extinction rather than a recovery later on). However, as a first approximation, the crude Lotka-Volterra model does give a sensible description. 

However, precisely because the population number oscillates, it is difficult to conclude whether the predator-prey scenario sufficiently reduces the number of total probes. It could be that we are near the ``valleys'' of the population number, when both predator and prey are relatively low in numbers. This is why the analysis in \cite{forgan} needs to be so thorough as to considering migration effect. Intuitively this amounts to taking the observer (us) into account: what is the probability that at the corner of the galaxy, we see very few SRPs, given the oscillatory populations? It is remarkable that the result turns out to be inadequate to explain Fermi paradox. 

In this work, we shall re-consider the application of the ``plain vanilla'' Lotka-Volterra equation. This requires some thoughts on the design of SRPs and how mutation might affect them. Assuming that errors are accumulated as the probes self-replicate (just like biological cell division, there might be programs in place to correct these errors, but errors are unavoidable), the new generations of SRPs would be expected to eventually be different from the original population in one way or another. In order to prevent the probes from consuming each other, the extraterrestrial engineers should design SRPs to still recognize mutated probes as self. Hence a probe can only become predatory when the self-recognition feature breaks down due to a new mutation and no longer recognizes other variants as self. 

Just like in \cite{forgan}, we will assume that the prey probes do not consume the predators. This assumption is not necessary. A predator is defined to be the population that benefits from the inter-species interaction and similarly the prey is the one that suffers. Therefore, the progenitor probe can prey on the predator probes, with some coefficient $+\beta' xy$ added into the Lotka-Volterra equation for $\d x/\d t$, but the overall coefficient $\beta'-\beta$ remains negative by definition of being a prey. Thus it is more convenient to just re-define $-\beta'+\beta \mapsto \beta$, or equivalently assuming $\beta'=0$.

With the above, the distinction between predator and prey, as in Eq.(\ref{vanilla}), can be made. In \cite{forgan}, the author mentioned the possibility of an ``omnivorous'' predator species that also consumes each other, and conclude that this should not affect the conclusion. Indeed, omnivorous predators are also \emph{replicating} while consuming, so this makes sense. Thus, we can assume without loss of generality that predator species would recognize \emph{themselves} and not become omnivorous. On the other hand, to replicate is the main purpose of SRPs. It is unlikely that the mutation into predator would also mutate the probes in so specific a way to make the predator probes \emph{exclusively} consuming the prey probes. (In fact, as we will see, extraterrestrial engineers should try their best to program SRPs to prevent exclusive predation behavior from arising.) In other words, \emph{the Lotka-Volterra equation should include also the resource term for the predators}. Unlike, say, foxes and rabbits in a typical ecosystem, in which the predators only eat the preys and would die off without the latter, SRPs dynamics would be quite different due to both species have access to environmental resources. This crucial term, missing in \cite{forgan}, makes all the difference.

A more realistic model should therefore be (we use Latin coefficients in place of Greek ones for the newly considered terms):
\begin{equation}
\left\{ \begin{aligned}
\frac{\d x}{\d t} &= \alpha x - \beta xy -Bx, ~~ \alpha > B,  \\
\frac{\d y}{\d t} &= \delta xy - \gamma y + Ay + Bx, ~~ A > \gamma,
\end{aligned} \right.
\end{equation}
where $A, B > 0$. Here $A$ is the coefficient of the resource term for the predator population. We expect that the magnitude of $\alpha$ is comparable to $A$, but they may not be identical after the mutation (the result below holds as long as $A > \gamma$). We have also included the mutation effect via the coefficient $B$, where preys are converted into predators with some fixed rate\footnote{This is akin to how the infected population is ``converted'' into the recovered population in a contagious illness ``SIR'' model \cite{Hethcote}, such as the Kermack-McKendrick model \cite{KM}. Note that mutation need not necessarily occurs \emph{during} replication, but also when other processes are taking place.}. The magnitude of $B$ should be much smaller than the other terms.
The condition $A > \gamma$ is simply the statement that the natural death rate is smaller than the growth rate of the population.

To solve for the fixed points, we solve the equations
\begin{flalign}
&x[(\alpha-B)-\beta y]=0 \notag \\
&y[\delta x + (A-\gamma)]+Bx=0.
\end{flalign}
This yields one trivial fixed point $(0,0)$, as well as one nontrivial fixed point: 
\begin{equation}
(x,y)=(x^*,y^*):=\left(\frac{(B-\alpha)(A-\gamma)}{B\beta + \delta(\alpha-B)},\frac{\alpha-B}{\beta}\right).
\end{equation} 
We note that the nontrivial fixed point is irrelevant because $x^* < 0$.

The Jacobian associated to the dynamical system is
\begin{equation}
\mathcal{J}=\begin{bmatrix} (\alpha-B)-\beta y & -\beta x \\ \delta y + B & \delta x -\gamma+A \end{bmatrix}.
\end{equation}

\begin{figure}
\begin{center}
\includegraphics[width=0.42\textwidth]{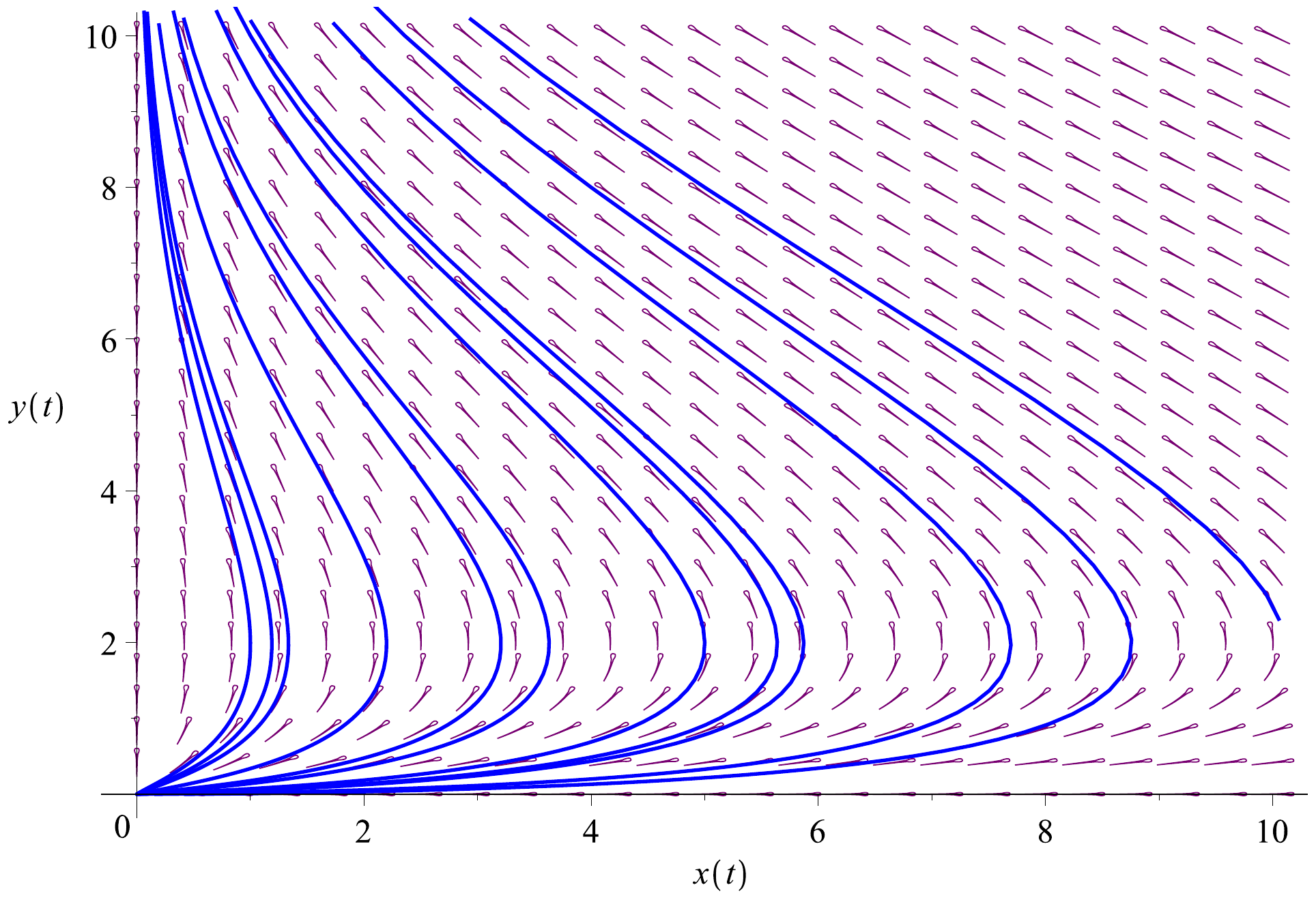}
\caption{The phase diagram of the improved predator-prey model. Here we use $\alpha=1,\beta=0.5,\gamma=0.2,\delta=0.1$, $B=0.01, A=0.95$. The origin -- the only fixed point -- is an unstable node (repeller). We observe that $x(t)\to 0$, i.e. the progenitor probes are driven into extinction.}
\label{fig3}
\end{center}
\end{figure}

\begin{figure}
\begin{center}
\includegraphics[width=0.40\textwidth]{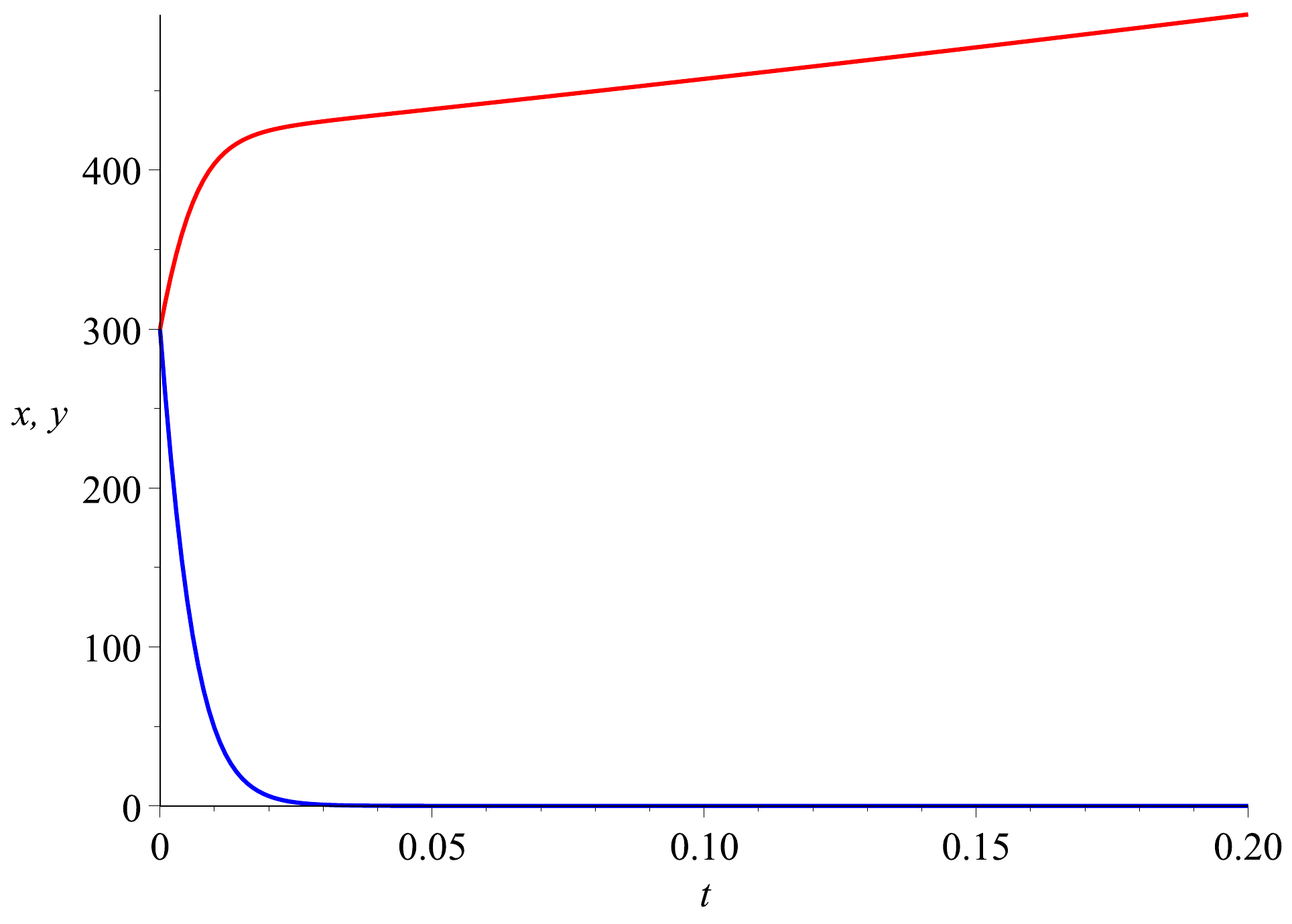}
\includegraphics[width=0.40\textwidth]{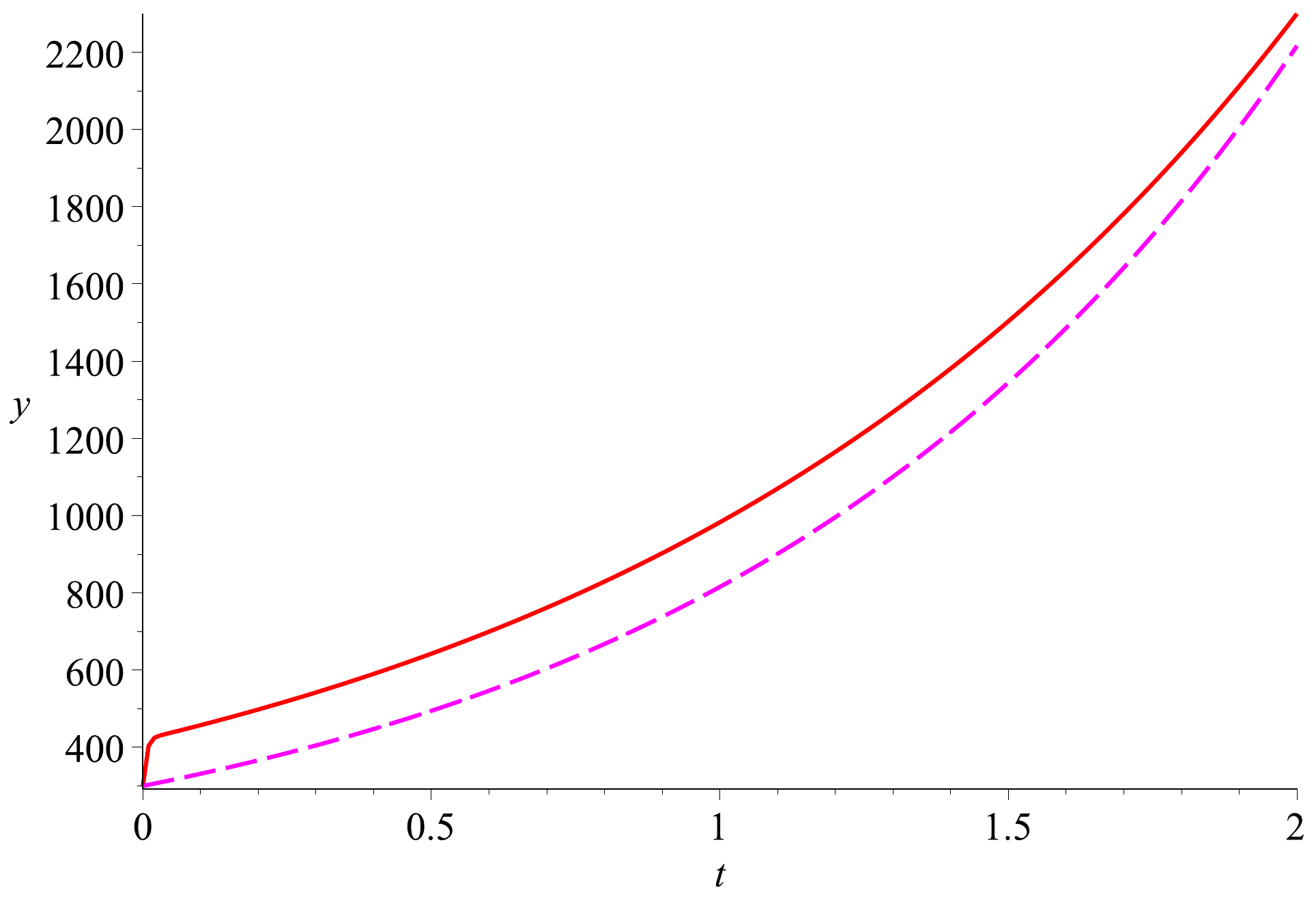}
\caption{\textbf{Top:} The number of the predator ($y$, the upper curve in red) and prey ($x$, the lower curve in blue) populations in the improved predator-prey model, with $\alpha=1,\beta=0.5,\gamma=0.2,\delta=0.1$, $B=0.01, A=0.95$. The initial conditions are $x(0)=y(0)=300$. \textbf{Bottom:} The subsequent exponential growth of $y(t)$ is comparable to the case of single-species probe without mutation (dashed).}
\label{fig4}
\end{center}
\end{figure}

The origin $(0,0)$ is the only fixed point, and since $\mathcal{J}|_{(0,0)}$ is a lower triangular matrix, its eigenvalues are $\alpha - B > 0$, and $A-\gamma > 0$. The trivial fixed point is therefore an unstable node that repels all nearby trajectories in the phase space (Fig.(\ref{fig3})). The result is that the prey quickly becomes extinct, while the predator population eventually grows exponentially (or tend to a constant if a finite capacity is assumed, \`a la the logistic model), as shown in the top plot of Fig.(\ref{fig4}). The predator growth rate is comparable to the progenitor probe in the scenario in which no mutation occurs; see the bottom comparison plots of Fig.(\ref{fig4}). We remark that the conversion term $Bx$ is not crucial, the major influence is due to the resource term. Even if we take $B=0$ in the model (which corresponds to the mutation being an one-off event), the repelling behavior persists.

\begin{figure}
\begin{center}
\includegraphics[width=0.43\textwidth]{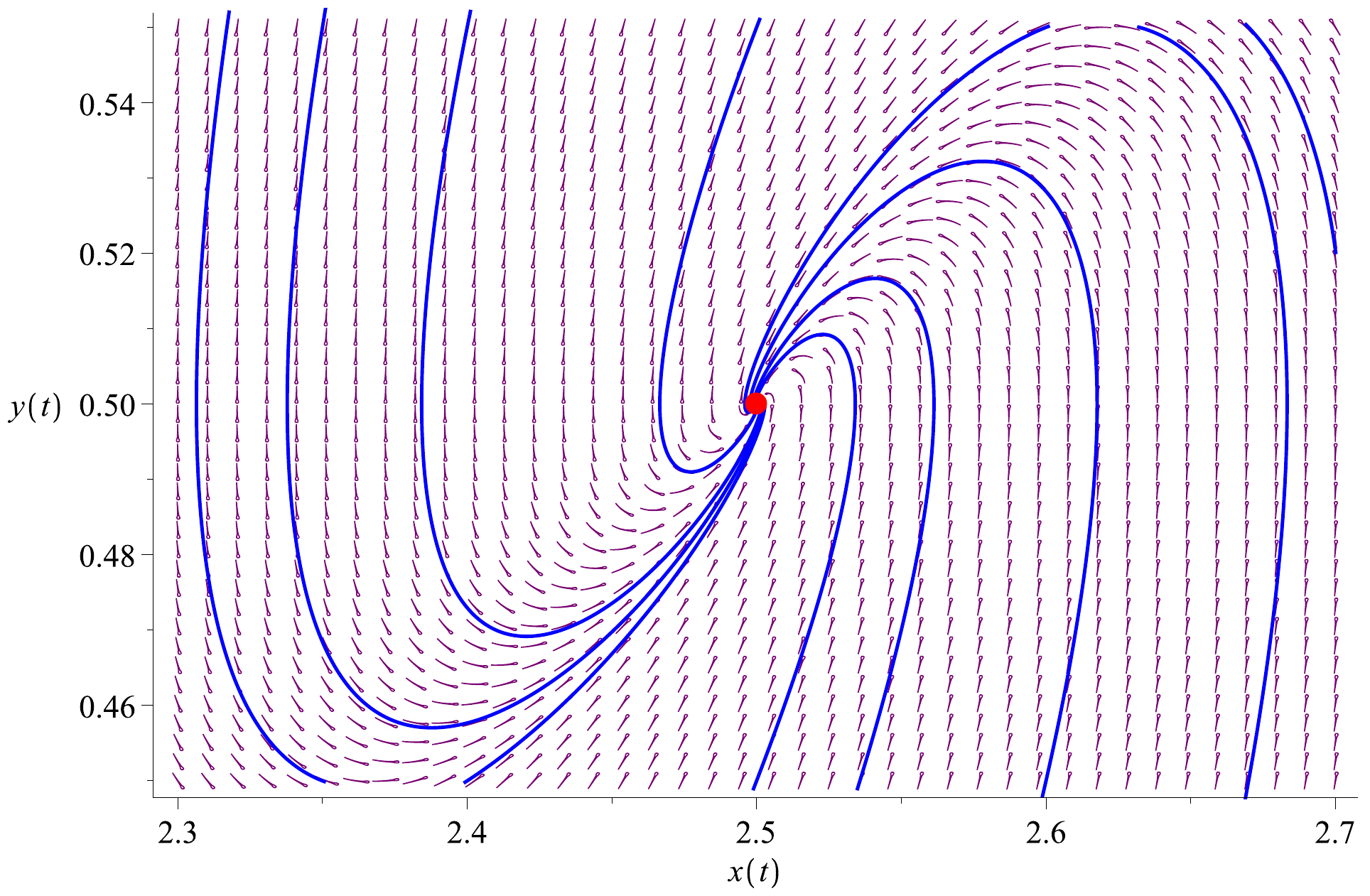}
\caption{The phase diagram of the predator-prey model with a ``defective'' predator. Here we use $\alpha=1,\beta=1,\gamma=5,\delta=1$, $B=0.5, A=0$. The system has a nontrivial fixed point marked with the red dot at $(x^*,y^*)=(2.5,0.5)$. In this case the trajectories are stable spirals, but this is not necessarily the case in general; see Fig.(\ref{fig6}).}
\label{fig5}
\end{center}
\end{figure}

\begin{figure}
\begin{center}
\includegraphics[width=0.43\textwidth]{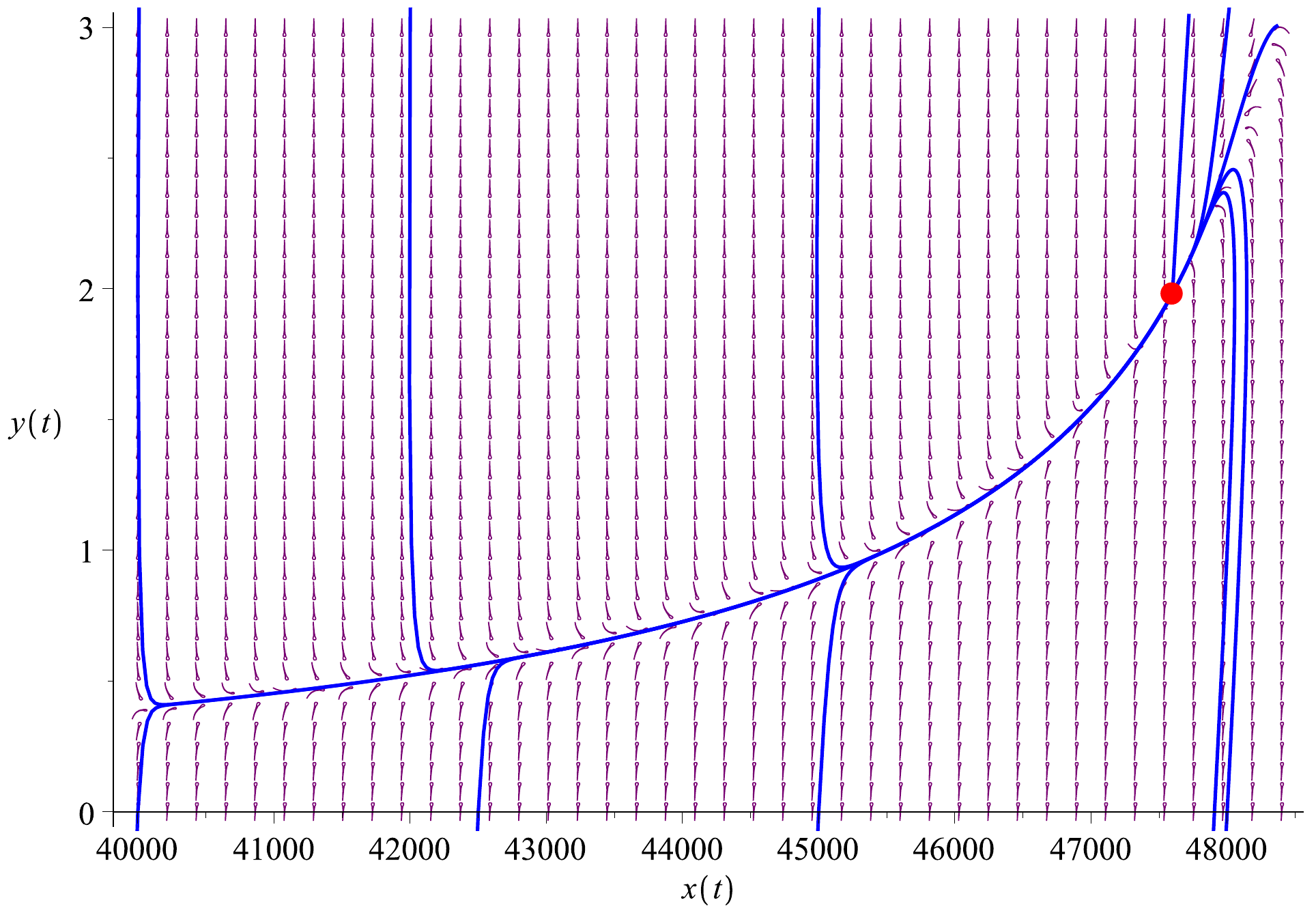}
\caption{The phase diagram of the predator-prey model with a ``defective'' predator. Here we use $\alpha=1,\beta=0.5,\gamma=5000,\delta=0.1$, $B=0.01, A=0$. The system has a nontrivial fixed point marked with a red dot at $(x^*,y^*)\approx (47596.15, 1.98)$, which is a stable node. The trajectories flow towards an attractor curve that leads to the stable node.}
\label{fig6}
\end{center}
\end{figure}

From our discussion above, it seems that we do not have to discuss the $A=0$ case (which is excluded with our assumption that $A > \gamma > 0$). However, this case will provide some insights into how to design SRPs in the Discussion section and thus is worth some investigations.
For $A=0$, the trivial fixed point $(0,0)$ is a saddle node. 
More crucially, the nontrivial fixed point is now relevant as $x^*,y^*>0$.
Though the system does not quite reduce to the ``plain vanilla'' Lotka-Volterra form (due to the $Bx$ term in the $\d y/ \d t$ equation), for small enough $B$, we still have oscillatory behaviors not too different from the  ``plain vanilla'' Lotka-Volterra equation, which requires a careful analysis of \cite{forgan}. This is not surprising since the term $Bx$ can be seen as a perturbation away from the  ``plain vanilla'' model. 
As $B$ gets larger, the nontrivial fixed point becomes an attractor. To see this, we note that the eigenvalues of $\mathcal{J}|_{(x^*,y^*)}(A=0)$ are
\begin{equation}
\lambda_{1,2}=-\frac{\gamma \beta B \pm \sqrt{\gamma^2 \beta^2 B^2 + \gamma f(\alpha,\beta,\delta,B)}}{2(B\beta + (\alpha - B)\delta)},
\end{equation}
where it can be shown that 
\begin{flalign}
 f(\alpha,\beta,\delta,B) = &-8B\beta\delta(B-\alpha)^2 + 4B^2\beta^2(B-\alpha)\notag \\ &+ 4\delta^2(B-\alpha)^3.
\end{flalign}
Each of the terms in $f(\alpha,\beta,\delta,B)$ are negative, so $f(\alpha,\beta,\delta,B)<0$. Consequently we see that indeed for $\gamma$ small enough, $\lambda_{1,2}$ are complex, whereas for large enough $\gamma$, both $\lambda_{1,2}<0$ so we have an attractor of the stable node type. 
For intermediate values of $\gamma$,  we would obtain a stable spiral attractor instead. 
An example is provided in Fig.(\ref{fig5}) with $\alpha=1,\beta=1,\gamma=5,\delta=1$, $B=0.5$, with eigenvalues $\lambda_{1,2} \approx -1.25 \pm 0.9682 i$. 
For the parameter values $\alpha=1,\beta=0.5,\delta=0.1, B=0.01$ that we use in Fig.(\ref{fig3}) and Fig.(\ref{fig4}), $\gamma$ needs to be of order $O(10^3)$ in order to exhibit an attractor behavior. An example with $\gamma=5000$ that gives rise to a stable node is given in Fig.(\ref{fig6}). Note that since the system is nonlinear, eigenvalue analysis for \emph{linear} stability need not quantitatively yield the correct result for when exactly the center transits into spirals.

For our purpose, we note that since $\gamma$ governs the predator population ``death rate'', this means that the attractor case corresponds to the predator probes being ``defective'', i.e., easily mulfunctioned.
This could happen as part of the programming. We will return to this issue in the Discussion. (The possibility that SRPs eventually mutate so much they cease to replicate was also discussed in \cite{1605.02169})
Note that the case of $A \neq 0$ but $A < \gamma$ is similar.

\section{Discussion: How (Not) to Build Replicators}
In this work we have shown that in a more realistic Lotka-Volterra model, a predatory SRP would drive its progenitor species towards extinction. In other words, the original SRPs would be replaced by the mutated version. Presumably this process will happen again and again since mutation would eventually give rise to new predators. The growth rate of the SRPs is therefore still exponential. This makes it clear that the Fermi paradox cannot be resolved.

This analysis also suggests a somewhat counter-intuitive way to program SRPs if we ever could use it to explore the Galaxy.
We might be tempted to think that we should design them to stop foraging as soon as their self-recognition program becomes sufficiently faulty. For example, ``self'' can be defined as a neighborhood in the parameter space, so that a small mutation is still acceptable. A ``kill'' code that renders the probe unable to forage or/and inactive can be designed to switch on when a newly constructed probe mutates outside the acceptable margin. This mimics a biological scenario: a cell which has an error in its DNA that cannot be repaired may undergo self-destruction (apoptosis), though this does not always work, which could lead to cancer. This parallels the runaways production of predator probes in our case. It is possible that after mutation the kill code is not as effective as to render the probe immediately inactive, but still the code could increase the ``death rate'' of the probe.

Such a design, however, can be risky. If consuming resources is coupled with self-recognition, the two behaviors are not independent, which conceivably increases the chance that a mutation might affect both. 
A biological analogy is the mutation of a virus -- when a virus mutates into a new variant that is more effective at spreading, its virulence may increase (or decrease) if this is coupled with transmissibility (which can either happen via a pleiotropic side effect or through genetic linkage), which is of course one current concern for COVID-19 \cite{covid}. 
For our case, a coupling between self-recognition and foraging behavior might not only \emph{cause} a predator behavior in the first place, but also conceivably makes the predator solely feeds on the progenitor probes $(A=0)$. In other words, such a design may increase the chance that an almost plain-vanilla predator-prey model emerges, at which point the growth rate of the probe population ceases to be exponential, and instead locked into oscillatory cycles, or tend to an attractor (if the predator is prey-specific \emph{and} defective in the sense that its ``death rate'' is high). The best strategy is not to program anything, so that foraging and self-recognition stay as independent systems. If predation behavior occurs by chance, it is less likely to be progenitor-specific. The predator probes would then out-compete the progenitors and grow exponentially -- which is the wanted behavior for SRPs. The risk of constructing SRPs is of course well-known and should be taken into account very carefully by any ambitious civilization \cite{safe}.

In any case, it does seem that based on our work and that of \cite{forgan}, mutated probes cannot explain the Fermi paradox. Instead, one is faced with the possibility that the lack of SRPs might indicate that highly advanced extraterrestrial civilizations do not exist \cite{ellery2}. Another possibility is that the probes are around, we just have not detected them \cite{1111.6131, freitas}. In fact, can we recognize a sufficiently advanced SRP technology if we see it? Of course, if SRPs do exist, their ``ecosystem'' might be very complicated depending on how they are programmed. Perhaps they are almost biological in the sense that mutation can lead to natural selection and vastly improved probes of many different species can evolve throughout the generations. One also has to take into account that there are more than one extraterrestrial civilizations that are launching SRPs, which would treat each other as resources to be consumed. We can even entertain the science fiction possibility that an extraterrestrial civilization would anticipate rival probes in advance and program their SRPs to be able to evade enemies and/or defend themselves; or more interestingly, establish contact and collaborate. 

Lastly, let us also mention the possibility of ``information panspermia'' first introduced by Gurzadyan \cite{pan}. This is the idea that information regarding the genomic sequence of an organism can be compressed by certain algorithm to be spread through the galaxy, either by transmitting (with the hope that some advanced civilization can download it) or by utilizing the SRPs, either by transmitting through the SRP network once it is established, or directly programmed into the SRPs. In the latter scenario, we can even envision that the probes from a sufficiently advanced civilization can ``seed'' life across the galaxy when they detected a suitable environment (hence conceivably re-creating even an ecosystem instead of just a single species), or even ``terraforming'' a targeted planet. Indeed, instead of fleets of spaceships, it is much more convenient to colonize the galaxy in this way. Though mutation to the SRPs would not affect its spread as is evident in the above analysis, it is possible that the mutation could also lead to certain corruptions of the genetic information. Depending on the degree of the corruption (and the efficiency of the error correction codes), this could either render the mutated probe unable to perform their panspermia function, or just cause the life seeded in this way to be genetic variations of the original species. If the latter scenario happens, it would be evolution in action on the galactic scale.

All these considerations, albeit some are highly speculative, imply that simple models based on Lotka-Volterra equations are unlikely to tell the whole story with regards to the Fermi paradox. Nevertheless, lacking any knowledge about the technology and the psychology of extraterrestrial species, we should always start from understanding the simplest models.

\begin{acknowledgments}
YCO thanks the National Natural Science Foundation of China (No.11922508) for funding support. YCO thanks Brett McInnes for useful discussions.
\end{acknowledgments}

\section*{Data Availability Statement}
Data sharing not applicable to this article as no datasets were generated or analysed during the current study.

\section*{Conflict of interest}
There is no conflict of interest.

\end{document}